\documentclass[journal]{IEEEtran}
\usepackage{amsmath,amssymb,amsfonts,amsthm,mathtools,cuted}
\usepackage{graphicx}
\usepackage[inkscapeversion=1.2.1]{svg}
\svgsetup{inkscapeexe=C:/Inkscape/bin/inkscape}
\usepackage{units}

\usepackage{import}
\usepackage{xspace}
\usepackage{pst-node,pst-plot}
\usepackage{mathtools}
\usepackage{caption}
\usepackage{subcaption}
\usepackage{xcolor}
\usepackage{xtab,booktabs}
\usepackage{longtable}
\usepackage{multicol}
\usepackage{dsfont}
\usepackage{bm}

\usepackage{algorithm}
\usepackage[noend]{algpseudocode}
\allowdisplaybreaks
\newcommand{\st}{\text{s.\,t. }}
\usepackage{hyperref}

\renewcommand{\vec}[1]{\mathbf{#1}} % Vectors
\newcommand{\gvec}[1]{\bm{#1}} % Greek vectors
\newcommand{\vbody}[2]{{\mathbb{#1}^{n_{#2}}}} % Vector body
\ifCLASSINFOpdf
\else
\fi

\begin{document}
\title{Federated $ K $-Means Clustering via Dual Decomposition-based Distributed Optimization}

\author{Vassilios~Yfantis,~Achim~Wagner,~Martin~Ruskowski % <-this % stops a space
\thanks{V. Yfantis and M. Ruskowski are with the Chair of Machine Tools and Control Sytems, Rheinland-Pfälzische Technische Universität Kaiserslautern-Landau (e-mail: vassilios.yfantis@mv.uni-kl.de).}% <-this % stops a space
\thanks{A. Wagner and M. Ruskowski are with the German Research Center for Artificial Intelligence.}% <-this % stops a space
\thanks{The authors would like to thank Vinit Hegiste for providing feedback on the manuscript.}
\thanks{The research was funded by a project supported by the Federal Ministry for Economic Affairs and Climate Action (BMWK) on the basis of a decision by the German Bundestag.}}

\maketitle

\begin{abstract}
The use of distributed optimization in machine learning can be motivated either by the resulting preservation of privacy or the increase in computational efficiency. On the one hand, training data might be stored across multiple devices. Training a global model within a network where each node only has access to its confidential data requires the use of distributed algorithms. Even if the data is not confidential, sharing it might be prohibitive due to bandwidth limitations. On the other hand, the ever-increasing amount of available data leads to large-scale machine learning problems. By splitting the training process across multiple nodes its efficiency can be significantly increased. This paper aims to demonstrate how dual decomposition can be applied for distributed training of $ K $-means clustering problems. After an overview of distributed and federated machine learning, the mixed-integer quadratically constrained programming-based formulation of the $ K $-means clustering training problem is presented. The training can be performed in a distributed manner by splitting the data across different nodes and linking these nodes through consensus constraints. Finally, the performance of the subgradient method, the bundle trust method, and the quasi-Newton dual ascent algorithm are evaluated on a set of benchmark problems. While the mixed-integer programming-based formulation of the clustering problems suffers from weak integer relaxations, the presented approach can potentially be used to enable an efficient solution in the future, both in a central and distributed setting.
\end{abstract}


% \fi
\IEEEpeerreviewmaketitle

\section{Introduction}
Training a machine learning model of any kind on a large set of data usually involves the solution of a challenging optimization problem. If the underlying data set becomes too large, it might not be possible to solve the resulting optimization problem in a reasonable amount of time. Distributed optimization methods can aid in rendering the optimization problem tractable through the use of multiple computational resources. Peteiro-Baral and Guijarro-Bedi\~{n}as provide an overview of methods for distributed machine learning \cite{Peteiro2013}. To train a global model in a distributed manner a consensus has to be established between the involved nodes and their underlying optimization problems. Forero et al. \cite{Forero2010,Forero2011} and Georgopoulos and Hasler \cite{Georgopoulos2014} demonstrate the distributed training of machine learning models using consensus-based distributed optimization. Tsianos et al. discuss practical issues with a consensus-based approach that arise from the difference between synchronous and asynchronous communication \cite{Tsianos2012}. Nedi\'{c} provides an overview of distributed gradient methods for convex training problems \cite{Nedic2020} while Verbraeken et al. give a general survey of distributed machine learning \cite{Verbraeken2020}.
\begin{figure}[h!]
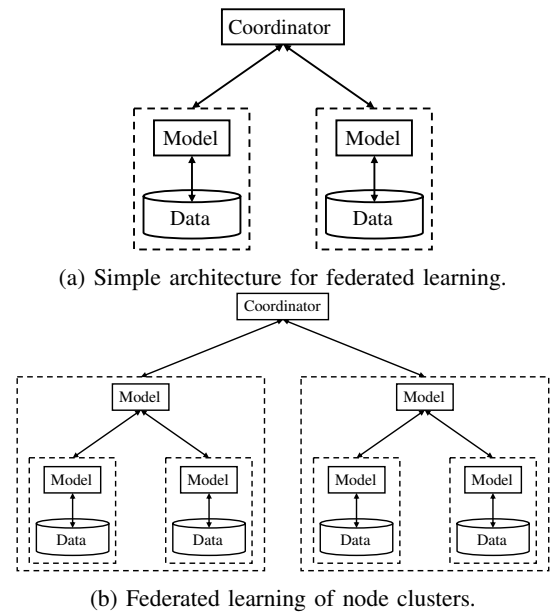

	\begin{subfigure}[b]{\linewidth}
		\centering
		\fontsize{8pt}{3pt}\selectfont
		\includesvg[width = 0.45\linewidth]{Figures/Clustering/FL_1}
		\caption{Simple architecture for federated learning.}
		\label{fig:fl_1}
	\end{subfigure}
	\begin{subfigure}[b]{\linewidth}
		\centering
		\fontsize{6pt}{3pt}\selectfont
		\includesvg[width = 0.8\linewidth]{Figures/Clustering/FL_2}
		\caption{Federated learning of node clusters.}
		\label{fig:fl_2}
	\end{subfigure}
	\label{fig:fl}
	\caption{Examples of federated learning architectures.}
\end{figure}

While computational performance remains an issue for many machine learning problems, the increase in computing power and in the efficiency of optimization algorithms can render many challenging problems tractable. However, the inability to share data due to confidentiality reasons still necessitates the use of distributed algorithms. Fig.~\ref{fig:fl_1} shows a setting in which training data is stored across two different nodes. Each node can use its local data to train an individual machine learning model. By including a coordination layer the two training processes can be guided in a way that a global model is trained, without the need to share confidential data. Distributed training of a global model without sharing individual training data is often referred to as federated optimization or federated learning \cite{Konevcny2016}. If the underlying optimization problems are still hard to solve, or if the training data is heterogeneous, the training process can be further divided into subproblems. Fig.~\ref{fig:fl_2} depicts the situation in which models of different node clusters are trained in a distributed manner which in turn are again coordinated to obtain a global model. This approach is often referred to as clustered federated learning \cite{Ghosh2020}.\\

Most algorithms for federated learning involve an averaging step of the model parameters of the individual nodes \cite{Mcmahan2017}. Federated learning methods have been applied in the context of manufacturing \cite{Hegiste2022}, healthcare \cite{Antunes2022}, mobile devices \cite{Lim2020}, and smart city sensing \cite{Jiang2020}. Li et al. \cite{Li2020} and Liu et al. \cite{Liu2022} provide surveys on federated learning while Chamikara et al. examine the privacy aspects related to external attacks \cite{Chamikara2021}. Applying federated learning to heterogeneous data sets can lead to the deterioration of the model quality of individual nodes in regard to their training data, which might hinder their willingness to participate in such a setting. This issue is addressed through personalized federated learning \cite{Kulkarni2020,Tan2022}.\\

$ K $-means clustering describes an unsupervised machine learning problem in which a set of observations/data is divided into $ K $ disjoint clusters according to a similarity measure \cite{Gambella2021}. Clustering problems can be found in many practical applications such as image segmentation \cite{Dhanachandra2015}, customer market segmentation \cite{Kansal2018}, or the identification of similar operating points in a production plant \cite{Rahimi2019}. While $ K $-means clustering is a well-studied problem, federated clustering, i.e., clustering of data across multiple nodes, has not been studied extensively in the literature yet. Dennis et al. present a one-shot federated clustering algorithm with heterogeneous data, where the clustering problem of each node is solved by Lloyd's heuristic algorithm \cite{Dennis2021heterogeneity, Lloyd1982}. Kumar et al. apply federated averaging to pre-trained models on separate devices and present an update strategy when new data points are added \cite{Kumar2020}. Li et al. address the security aspect of federated clustering by encoding the data of each node and applying Lloyd's algorithm to the encoded data \cite{Li2022}. Stallmann and Wilbik extend fuzzy $c$-means clustering, i.e., a clustering problem where each data point can be assigned to multiple clusters, to a federated setting \cite{Stallmann2022}. A similar approach to federated $c$-means clustering was previously presented by Pedrycz \cite{Pedrycz2021}. Wang et al. use model averaging and gradient sharing for federated clustering of data in a smart grid \cite{Wang2022federated}.\\

A common feature of the federated clustering approaches described in the literature is the use of a heuristic algorithm to solve the individual clustering problems. In this paper mixed-integer programming is used to solve the individual clustering problems and a dual decomposition-based distributed optimization approach is employed to coordinate the solutions of the different nodes. While an averaging step is still performed to obtain feasible primal solutions, the use of duality enables the computation of valid lower bounds on the objective of the global clustering problem. It should be noted that the mixed-integer programming problems resulting from the $K$-means clustering problem are very hard to solve due to their weak integer relaxations. Thus, the approach presented in this paper can be regarded as preliminary work that can become a viable option with the continuous improvement of mixed-integer programming solvers \cite{Achterberg2013, Koch2022}.
\section{Notation}
We denote vectors with bold lower-case letters ($\vec{x}$) and matrices with bold upper-case letters ($\vec{X}$). The $j$th element of a vector $\vec{x}$ is denoted by $[\vec{x}]_j$. The Euclidean norm is denoted by $\Vert \cdot \Vert_2$. The value of a variable $\vec{x}$ in iteration $t$ is denoted by $\vec{x}^{(t)}$.

\section{$ K $-Means clustering}
 This section presents the mixed-integer programming-based formulation of the $ K $-means clustering training problem. The formulation is subsequently extended to the case of distributedly stored data, which gives rise to a federated learning problem. Consensus constraints are used to couple the training problems of different nodes. These constraints can be dualized such that the federated learning problem can be solved via dual decomposition-based distributed optimization. Since the underlying optimization problem contains integrality constraints it is not convex and thus strong duality does not hold. However, a feasible primal solution can be computed in each iteration through an averaging heuristic.
\subsection{MIQCP formulation}
\begin{figure*}[h!]
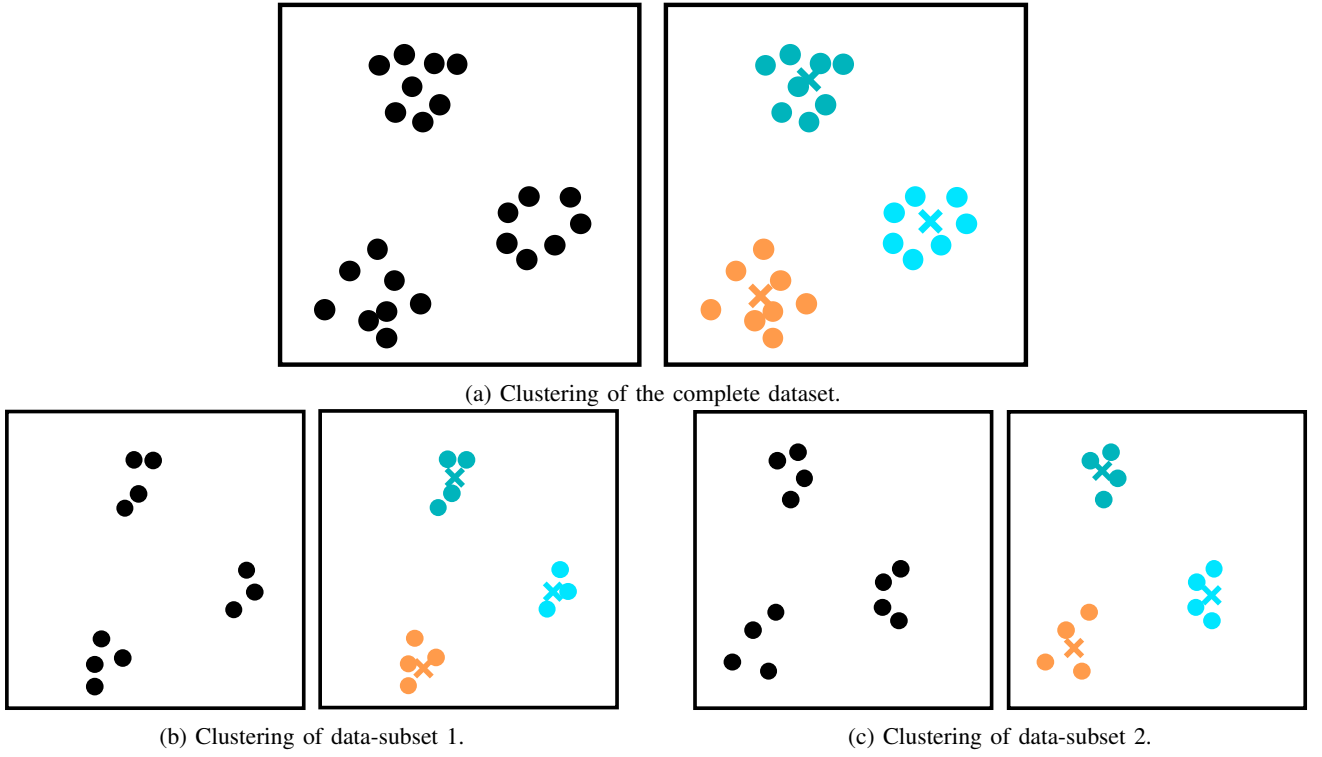

	\begin{subfigure}[t]{\linewidth}
		\centering
		\fontsize{14pt}{3pt}\selectfont
		\includesvg[width = 0.55\linewidth]{Figures/Clustering/Clustered_Data_0}
		\caption{Clustering of the complete dataset.}
		\label{fig:clustered_data_0}
	\end{subfigure}
	
	\begin{subfigure}[b]{0.5\linewidth}
		\centering
		\fontsize{10pt}{3pt}\selectfont
		\includesvg[width = 0.9\linewidth]{Figures/Clustering/Clustered_Data_1}
		\caption{Clustering of data-subset 1.}
		\label{fig:clustered_data_1}
	\end{subfigure}
	\begin{subfigure}[b]{0.5\linewidth}
		\centering
		\fontsize{10pt}{3pt}\selectfont
		\includesvg[width = 0.9\linewidth]{Figures/Clustering/Clustered_Data_2}
		\caption{Clustering of data-subset 2.}
		\label{fig:clustered_data_2}
	\end{subfigure}
	\label{fig:clustered_data}
	\caption{Illustration of $ K $-means clustering both in a centralized and a decentralized setting.}
\end{figure*}
The goal of $ K $-means clustering is to assign a set of observations $ \vec{y}_j \in\vbody{R}{\vec{y}},\;j \in \mathcal{J} = \{1,\dots,\vert\mathcal{J}\vert\} $ to a set of clusters $ \mathcal{K} = \{1,\dots,K\} $ and to compute the centroids of each cluster. The number of clusters is a hyper-parameter and is set a priori or in an iterative manner. This problem can be formulated as a mixed-integer nonlinear programming (MINLP) problem \cite{Aloise2012,Gambella2021},
\begin{subequations}\label{Clustering_MINLP}
	\begin{align}
		\underset{w_{jk},\vec{m}_k}{\min}&\sum_{j\in\mathcal{J}}\sum_{k\in\mathcal{K}} w_{jk}\cdot \Vert \vec{y}_j - \vec{m}_k \Vert_2^2,\\
		\st &\sum_{k\in\mathcal{K}}w_{jk} = 1,\forall j\in\mathcal{J},\\
		&w_{jk} \in \{0,1\},\;\forall j\in\mathcal{J},k\in\mathcal{K},\;\vec{m}_k\in\vbody{R}{\vec{y}}\;\forall k\in\mathcal{K}.
	\end{align}
\end{subequations}
The binary variables $ w_{jk} $ indicate if observation $ \vec{y}_j $ is assigned to cluster $ k $ and $ \vec{m}_k $ is the centroid of cluster $ k $. Constraint (\ref{Clustering_MINLP}b) enforces that each observation is assigned to exactly one cluster, while the objective is to minimize the sum of the squared Euclidean distances of all observations to the centroids of their assigned clusters.
Problem \eqref{Clustering_MINLP} is a nonconvex MINLP which is hard to solve. In practice it is more efficient to use a linearized formulation by introducing the variable $ d_{jk} $, which describes the squared distance between an observation $ j $ and the centroid of cluster $ k $ \cite{Gambella2021},
\pagebreak
\begin{subequations}\label{Clustering_MILP}
	\begin{align}
		\underset{w_{jk},d_{jk},\vec{m}_k}{\min}&\sum_{j\in\mathcal{J}}\sum_{k\in\mathcal{K}}d_{jk}\\
		\st &\sum_{k\in\mathcal{K}}w_{jk} = 1,\forall j\in\mathcal{J},\\
        \nonumber
		&d_{jk} \ge \Vert \vec{y}_j - \vec{m}_k \Vert_2^2 - M_j\cdot(1-w_{jk}),\\
        &\forall j\in\mathcal{J},k\in\mathcal{K}\\
        \nonumber
		&w_{jk} \in \{0,1\}, d_{jk}\ge0,\\
        &\forall j\in\mathcal{J},k\in\mathcal{K},\;\vec{m}_k\in\vbody{R}{\vec{y}}\;\forall k\in\mathcal{K}.
	\end{align}
\end{subequations}
\eqref{Clustering_MILP} is a mixed-integer quadratically constrained programming (MIQCP) problem with a convex integer relaxation. Constraint (\ref{Clustering_MILP}c) is an epigraph formulation of the squared Euclidean distance if observation $ j $ is assigned to cluster $ k $, i.e., when $ w_{jk} = 1 $. Otherwise, the parameter $ M_j $ has to be large enough so that the constraint is trivially satisfied for $ w_{jk} = 0 $. In theory, a common big-M parameter can be used for all constraints described by (\ref{Clustering_MILP}c). However, the parameter should be chosen as small as possible to avoid weak integer relaxations. In the following, the big-M parameter is set as
\begin{subequations}\label{big_M}
	\begin{alignat}{2}
    M_j = &\underset{\gvec{\chi} \in \mathcal{Y}}{\max}\;\Vert \vec{y}_j - \gvec{\chi}\Vert_2^2,\; \forall j \in \mathcal{J},\\
		\mathcal{Y} = &\{\vec{y}\in \vbody{R}{\vec{y}} \vert\; \underset{j \in \mathcal{J}}{\min}\; [\vec{y}_j]_l \le [\vec{y}]_l \le  \underset{j \in \mathcal{J}}{\max}\; [\vec{y}_j]_l\\
        &l=1.\dots,n_{\vec{y}}\}.
	\end{alignat}
\end{subequations}

Different approaches have been proposed to solve the clustering optimization problem. Bagirov and Yearwood present a heuristic method based on nonsmooth optimization \cite{Bagirov2006}, Aloise et al. propose a column generation algorithm \cite{Aloise2012} and Karmitsa et al. use a diagonal bundle method \cite{Karmitsa2017}. Fig.~\ref{fig:clustered_data_0} illustrates the concept of $ K $-means clustering. The unlabeled data (left) is split into 3 clusters according to their distance to the computed cluster centroid (crosses).
\subsection{Distributed consensus formulation}
Problem \eqref{Clustering_MILP} describes the case in which the entire data set is accessible from a single node. However, this might not always be the case, especially if the underlying data is confidential. In the following it is assumed that the data set is split across several nodes $ \mathcal{I} = \{1,\dots,N_s\} $, with each node having access to the data-subset $ \mathcal{J}_i \subset \mathcal{J} $. The MIQCP problem \eqref{Clustering_MILP} can be extended to the case of multiple nodes,
\begin{subequations}\label{Clustering_MILP_Central}
	\begin{align}
		\underset{w_{ijk},d_{ijk},\vec{m}_{k}}{\min}&\sum_{i\in\mathcal{I}}\sum_{j\in\mathcal{J}_i}\sum_{k\in\mathcal{K}}d_{ijk}\\
		\st &\sum_{k\in\mathcal{K}}w_{ijk} = 1,\forall i \in\mathcal{I}, j\in\mathcal{J}_i,\\
		\nonumber
        &d_{ijk} \ge \Vert \vec{y}_j - \vec{m}_{k} \Vert_2^2 - M_j\cdot(1-w_{ijk}),\\
        &\forall i \in\mathcal{I}, j\in\mathcal{J}_i,k\in\mathcal{K},\\
	       \nonumber	
        &w_{ijk} \in \{0,1\}, d_{ijk}\ge0,\;\forall i \in \mathcal{I}, j\in\mathcal{J}_i,k\in\mathcal{K},\\
        &\vec{m}_{k}\in\vbody{R}{\vec{y}},\;\forall k\in\mathcal{K}.
	\end{align}
\end{subequations}
The goal of problem \eqref{Clustering_MILP_Central} is again to compute a set of cluster centroids $ \vec{m}_k $ and to assign the observations of all nodes to these clusters. However, if the nodes cannot share their data, problem \eqref{Clustering_MILP_Central} cannot be solved in a centralized manner. A simple distributed approach would be to solve a clustering problem in each node $ i $. This could lead to a situation as depicted in Fig.~\ref{fig:clustered_data_1} and Fig.~\ref{fig:clustered_data_2}. If each the data set is split across two nodes, each one can solve a clustering problem. However, both nodes will compute different cluster centroids.\\
The goal of a federated learning approach is to train a global model, i.e., global cluster centroids in the case of $ K $-means clustering, without sharing the local data between the nodes. To this end each node $ i $ can compute individual cluster centroids $ \vec{m}_{ik} $,
\begin{subequations}\label{Clustering_MILP_Consenus}
	\begin{align}
		\underset{w_{ijk},d_{ijk},\vec{m}_{ik}}{\min}&\sum_{i\in\mathcal{I}}\sum_{j\in\mathcal{J}_i}\sum_{k\in\mathcal{K}}d_{ijk}\\
		\st &\sum_{k\in\mathcal{K}}w_{ijk} = 1,\forall i \in\mathcal{I}, j\in\mathcal{J}_i,\\
        \nonumber
		&d_{ijk} \ge \Vert \vec{y}_j - \vec{m}_{ik} \Vert_2^2 - M_j\cdot(1-w_{ijk}),\\
        &\forall i \in\mathcal{I}, j\in\mathcal{J}_i,k\in\mathcal{K},\\
		&\vec{m}_{ik} = \vec{m}_{i'k},\;\forall i \in \mathcal{I}, i' \in \mathcal{N}_i, k \in \mathcal{K},\\
        \nonumber
		&w_{ijk} \in \{0,1\}, d_{ijk}\ge0,\\
        &\forall i \in \mathcal{I}, j\in\mathcal{J}_i,k\in\mathcal{K},\\
        &\vec{m}_{ik}\in\vbody{R}{\vec{y}}\;\forall, i \in \mathcal{I}, k\in\mathcal{K}.
	\end{align}
\end{subequations}
Since the goal is to obtain global cluster centroids, the individual cluster centroids are coupled through consensus constraints (\ref{Clustering_MILP_Consenus}d), where $ \mathcal{N}_i $ contains the set of neighboring nodes of node $ i $. Problem \eqref{Clustering_MILP_Consenus} describes a set of $ N_s $ subproblems coupled through the consensus constraints. In the following subsection dual variables are used to decouple the clustering problems of the different nodes.
\section{Dual decomposition-based distributed clustering}
\begin{figure}[h!]
	\centering
	\fontsize{14pt}{3pt}\selectfont
	\includesvg[width = 0.7\linewidth]{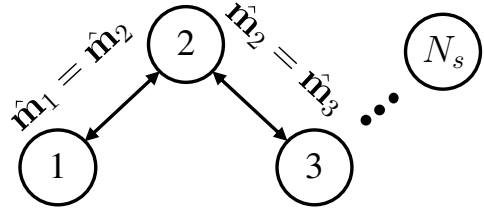}
	\caption{Illustration of a linear network topology and the resulting consensus constraints.} 
	\label{fig:cluster_network}
\end{figure}

This section presents how the consensus formulation (\ref{Clustering_MILP_Consenus}) of the clustering problem can be decomposed by introducing dual variables. Dual decomposition can be applied to constraint-coupled optimization problems of the form
\begin{subequations}\label{Constraint_Coupled}
	\begin{align}
		\underset{\vec{x}}{\min}\;&\sum_{i\in\mathcal{I}}f_i(\vec{x}_i),\\
		\st & \sum_{i\in\mathcal{I}}\vec{A}_i\vec{x}_i = \vec{b},\\
		&\vec{x}_i \in \mathcal{X}.
	\end{align}
\end{subequations}
Equation (\ref{Constraint_Coupled}) describes an optimization problem consisting of a set of $ {\mathcal{I}=\{1,\dots,N_s\}} $ subproblems. The subproblems are coupled through the constraints (\ref{Constraint_Coupled}b) and each one is described by individual variables $ \vec{x}_i $ and constraints $ \mathcal{X}_i $. Dual decomposition is based on the introduction of dual variables for the coupling constraints (\ref{Constraint_Coupled}b) and the solution of the resulting dual optimization problem. The idea was first introduced by Everett \cite{Everett.1963} for problems involving shared limited resources. Problem \eqref{Clustering_MILP_Consenus} can also be rewritten as a general constraint-coupled optimization problem by defining the matrix $ \vec{A} $ describing the connections between the different nodes. In the following only linear network topologies as depicted in Fig.~\ref{fig:cluster_network} are considered. Note that the discussion in the remainder of this paper can be easily extended to different network topologies.\\

By defining the vector of stacked cluster centroids of each node $ i $,
\begin{equation}\label{Stacked_Centroids}
	\hat{\vec{m}}_{i} \coloneqq \begin{bmatrix}
		\vec{m}_{i,1}\\
		\vdots\\
		\vec{m}_{i,k}
	\end{bmatrix} \in \mathbb{R}^{K\cdot n_{\vec{y}}},
\end{equation}
the consensus constraints can be rewritten as
\begin{subequations}\label{Consensus_Alg}
	\begin{align}
		&\hat{\vec{m}}_{1} - \hat{\vec{m}}_{2} = \vec{0},\\
		&\hat{\vec{m}}_{2} - \hat{\vec{m}}_{3} = \vec{0},\\
		\nonumber										  &\vdots\\								&\hat{\vec{m}}_{N_s -1} - \hat{\vec{m}}_{N_s} = \vec{0}.			  
	\end{align}
\end{subequations}
Constraints \eqref{Consensus_Alg} can subsequently be rewritten in matrix form
\begin{subequations}\label{Consensus_Matrix}
	\begin{align}
		\underbrace{\begin{bmatrix}
				\vec{I} & -\vec{I} & \vec{0} & \cdots & \vec{0} & \vec{0} \\
				\vec{0} & \vec{I} & -\vec{I} & \cdots & \vec{0} & \vec{0} \\
				\vdots & \vdots & \vdots & \ddots & \vdots & \vdots \\
				\vec{0} & \vec{0} & \vec{0} & \cdots & \vec{I} & -\vec{I}
		\end{bmatrix}}_{=:\vec{A} \in \mathbb{R}^{K\cdot n_{\vec{y}}\cdot(N_s-1) \times K \cdot n_{\vec{y}}\cdot N_s  } } \cdot \begin{bmatrix}
			\hat{\vec{m}}_1 \\
			\hat{\vec{m}}_2 \\
			\hat{\vec{m}}_3 \\
			\vdots \\
			\hat{\vec{m}}_{N_s}
		\end{bmatrix} = \vec{0},
	\end{align}
\end{subequations}
or in a more compact way
\begin{equation}\label{Consensus_Matrix_Compact}
	\sum_{i\in\mathcal{I}}\vec{A}_i\hat{\vec{m}}_{i} = \vec{0}
\end{equation}
with $ \vec{A}_i \in  \mathbb{R}^{K\cdot n_{\vec{y}}\cdot(N_s-1) \times K\cdot n_{\vec{y}} }$. By introducing dual variables $ \gvec{\lambda} \in \vbody{R}{K\cdot n_{\vec{y}}\cdot(N_s-1)} $  for the consensus constraints \eqref{Consensus_Matrix_Compact} the Lagrange function of problem (\ref{Clustering_MILP_Consenus}) can be defined,
\begin{equation}\label{Clustering_Lagrange}
	\mathcal{L}(w_{ijk},d_{ijk},\vec{m}_{ik},\gvec{\lambda}) =  \sum_{i\in\mathcal{I}}\underbrace{\left(\sum_{j\in\mathcal{J}_i}\sum_{k\in\mathcal{K}}d_{ijk} + \gvec{\lambda}^T \vec{A}_i\hat{\vec{m}}_{i}\right)}_{=:\mathcal{L}_i(w_{ijk},d_{ijk},\vec{m}_{ik},\gvec{\lambda})}.
\end{equation}
The minimization of the Lagrange function for a fixed value of the dual variables $ \gvec{\lambda} $ gives the corresponding value of the dual function.
	\begin{align}\label{Clustering_Dual}
		d(\gvec{\lambda}) \coloneqq 	\underset{w_{ijk},d_{ijk},\vec{m}_{ik}}{\min}&\sum_{i\in\mathcal{I}}\mathcal{L}_i(w_{ijk},d_{ijk},\vec{m}_{ik},\gvec{\lambda})\\
  \nonumber
		\st & \text{(\ref{Clustering_MILP_Consenus}b), (\ref{Clustering_MILP_Consenus}c), (\ref{Clustering_MILP_Consenus}e), (\ref{Clustering_MILP_Consenus}f)}
	\end{align}
The dual function has two important properties. First, the value of the dual function is always a lower bound on the solution of its corresponding primal problem, in this case, problem (\ref{Clustering_MILP_Consenus}) \cite{Nocedal2006}. The problem of finding the dual variables that result in the best lower bound is referred to as the dual optimization problem,
\begin{equation}\label{Dual_Problem}
	\underset{\gvec{\lambda}}{\max}\;d(\gvec{\lambda}).
\end{equation}
The resulting dual problem can be solved in a distributed manner by solving the individual clustering problems for the current value of the dual variables,
\begin{subequations}\label{Individual_Clustering_Problem}
	\begin{align}
		\underset{w_{ijk},d_{ijk},\vec{m}_{ik}}{\min}&\mathcal{L}_i(w_{ijk},d_{ijk},\vec{m}_{ik},\gvec{\lambda})\\
		\st &\sum_{k\in\mathcal{K}}w_{ijk} = 1,\;\forall j\in\mathcal{J}_i,\\
        \nonumber
		&d_{ijk} \ge \Vert \vec{y}_j - \vec{m}_{ik} \Vert_2^2 - M_j\cdot(1-w_{ijk}),\\
        &\forall j\in\mathcal{J}_i,k\in\mathcal{K},\\
		\nonumber
		&w_{ijk} \in \{0,1\}, d_{ijk}\ge0,\;\forall  j\in\mathcal{J}_i,k\in\mathcal{K},\\
		&\vec{m}_{ik}\in\vbody{R}{\vec{y}}\;\forall k\in\mathcal{K}.
	\end{align}
\end{subequations}
Second, the dual function (\ref{Clustering_Dual}) is always concave, regardless of whether the primal problem is convex or not \cite{Nocedal2006}. Therefore the dual problem (\ref{Dual_Problem}) is a convex optimization problem. However, the dual function is usually nondifferentiable due to a changing set of active individual constraints, which means that problem (\ref{Dual_Problem}) is a nonsmooth optimization problem \cite{Yfantis2023EURO}. The following subsections present some algorithms for the solution of the dual problem, namely the subgradient method, the bundle trust method, and the quasi-Newton dual ascent algorithm.
\subsection{Subgradient method}
Since the dual function is nondifferentiable a gradient cannot be defined for every value of the dual variables. Instead, a subgradient can be used. A vector $ \gvec{\xi}\in\vbody{R}{\gvec{\chi}} $ is a subgradient of a concave function $ \phi(\gvec{\chi}) $ at a point $ \gvec{\chi}_0 $ if
\begin{equation}\label{Subgradient_Condition}
	\phi(\gvec{\chi}) \le \phi(\gvec{\chi}_0) + \gvec{\xi}^T(\gvec{\chi}-\gvec{\chi}) 
\end{equation}
 for all $ \gvec{\chi}\in \text{dom}\;\phi $. The set of all subgradients at a point $ \gvec{\chi}_0 $ comprise the subdifferential $ \partial \phi(\gvec{\chi}_0) $ Technically equation (\ref{Subgradient_Condition}) defines a supergradient. Nevertheless, the term subgradient is commonly used in the literature for both convex and concave functions.\\
 A subgradient of the dual function for a given value of the dual variables $ \gvec{\lambda}^{(t)} $ can be computed by evaluating the coupling constraints (\ref{Consensus_Matrix_Compact}),
 \begin{equation}\label{Subgradient_Consensus}
 	\vec{g}(\gvec{\lambda}^{(t)}) = \sum_{i\in\mathcal{I}}\vec{A}_i\hat{\vec{m}}_{i}(\gvec{\lambda}^{(t)}) \in \partial d(\gvec{\lambda}^{(t)}),
 \end{equation}
where $ \hat{\vec{m}}_{i}(\gvec{\lambda}^{(t)}) $ are the cluster centroids obtained by solving the individual clustering problems (\ref{Individual_Clustering_Problem}).\\
In the subgradient method the dual variables are updated in each iteration $ t $ along the direction of the subgradient \cite{Shor2012minimization}
\begin{equation}\label{Subgradient_Update}
	\gvec{\lambda}^{(t+1)} = \gvec{\lambda}^{(t)} + \alpha^{(t)}\vec{g}(\gvec{\lambda}^{(t)}),
\end{equation} 
where $ \alpha^{(t)} $ is a step size parameter. The step size parameter plays an important role in the convergence of the algorithm. If it is chosen too large the algorithm might diverge, while a too small choice might significantly slow down its convergence. A common choice to adapt the step size throughout the iterations is
\begin{equation}\label{Step_Size}
	\alpha^{(t)} = \alpha^{(0)}/\sqrt{t},
\end{equation}
with an initial step size $ \alpha^{(0)} $ \cite{Bertsekas.1999}.
\subsection{Bundle trust method}
The subgradient method usually exhibits a slow rate of convergence, since only using information from the current subgradient may not provide an ascent direction for the algorithm. Bundle methods are generally more efficient by utilizing multiple subgradients from previous iterations \cite{Makela.2002}. To this end the data
\begin{multline}\label{Bundle}
	\mathcal{B}^{(t)} = \{(\gvec{\lambda}^{(l)}, \vec{g}(\gvec{\lambda}^{(l)}), d(\gvec{\lambda}^{(l)})) \in \vbody{R}{\gvec{\lambda}}\times \vbody{R}{\gvec{\lambda}} \times \mathbb{R}\\
    \vert\;l=t-\tau+1,\dots, t	\}
\end{multline}
is stored in each iteration, where $ n_{\gvec{\lambda}} $ denotes the number of dual variables. $ \mathcal{B}^{(t)} $ is referred to as a bundle and it contains the dual variables, subgradients, and values of the dual function from previous iterations. Since storing all information from all previous iterations might cause memory issues, only data from the previous $ \tau $ iterations is used.\\
The idea of bundle methods is to use the collected information to construct a piece-wise linear over-approximation of the nonsmooth dual function $ d(\gvec{\lambda}) $, a so-called cutting plane model,
\begin{equation}\label{CPM}
	\hat{d}^{(t)}(\gvec{\lambda}) \coloneqq \underset{l\in \{t-\tau+1,\dots, t\}}{\min}\{d(\gvec{\lambda}^{(l)}) + \vec{g}^T(\gvec{\lambda}^{(l)})(\gvec{\lambda} - \gvec{\lambda}^{(l)})\}.
\end{equation}
The approximation can be written in an equivalent form as
\begin{equation}\label{CPM_LinError}
	\hat{d}^{(t)}(\gvec{\lambda}) = \underset{l\in \{t-\tau+1,\dots, t\}}{\min}\{d(\gvec{\lambda}^{(t)}) + \vec{g}^T(\gvec{\lambda}^{(l)})(\gvec{\lambda} - \gvec{\lambda}^{(t)}) - \beta^{(l,t)}\},
\end{equation}
with the linearization error
\begin{multline}\label{LinError}
	\beta^{(l,t)} = d(\gvec{\lambda}^{(t)}) - d(\gvec{\lambda}^{(l)})- \vec{g}^T(\gvec{\lambda}^{(l)})(\gvec{\lambda}^{(t)} - \gvec{\lambda}^{(l)}),\\
 \forall l\in \{t-\tau+1,\dots, t\}.
\end{multline}
The update direction of the dual variables can then be computed by solving a direction-finding problem
\begin{subequations}\label{NonsmoothDirectipnFinding}
	\begin{align}
		\underset{\vec{s
			}\in \vbody{R}{\gvec{\lambda}}}{\max}\;&\hat{d}^{(t)}(\gvec{\lambda}^{(t)} + \vec{s}),\\
		\text{s.\ t. } & \Vert \vec{s} \Vert_2^2 \le \alpha^{(t)},\\
	\end{align}
\end{subequations}
where constraint (\ref{NonsmoothDirectipnFinding}b) represents a trust region. Therefore, this variant of the bundle method is referred to as the bundle trust method (BTM). Other variants include proximal bundle methods, where the trust region is replaced by a regularization term in the objective function \cite{Bagirov2014}. Problem (\ref{NonsmoothDirectipnFinding}) is still a nonsmooth optimization problem and can be transformed into a smooth quadratic direction finding problem by using an epigraph formulation,\pagebreak
\begin{subequations}\label{SmoothDirectipnFinding}
	\begin{align}
		\underset{v\in\mathbb{R},\;\vec{s} \in \vbody{R}{\gvec{\lambda}}}{\max}\;&v,\\
		\text{s.\ t. } & \Vert \vec{s} \Vert_2^2 \le \alpha^{(t)},\\
		&\vec{g}^T(\gvec{\lambda}^{(l)})\vec{s} - \beta^{(l,t)} \ge v,\;\forall l\in \{t-\tau+1,\dots, t\}.
	\end{align}
\end{subequations}
After computing a direction the dual variables are updated according to
\begin{equation}\label{BTM_Update}
	\gvec{\lambda}^{(t+1)} = \gvec{\lambda}^{(t)} + \vec{s}^{(t)}.
\end{equation}
Bundle methods are widely used in machine learning, as nonsmoothness is encountered in many training problems involving regularization terms \cite{Le2007}. Bundle methods can also be used to solve the clustering problem (\ref{Clustering_MILP}) \cite{Karmitsa2017}. However, note that in this paper the BTM algorithm is used to solve the nonsmooth dual problem (\ref{Dual_Problem}).
\subsection{Quasi-Newton dual ascent}
Since the dual function is always concave it can be locally approximated by a quadratic function. Yfantis et al.  recently proposed the quasi-Newton dual ascent (QNDA) algorithm that approximates the dual function by a quadratic function \cite{Yfantis2022, Yfantis2023EURO},
\begin{multline}\label{QNDA_Dual}
	d_B^{(t)}(\gvec{\lambda})=\frac{1}{2}(\gvec{\lambda}-\gvec{\lambda}^{(t)})^T \vec{B}^{(t)}(\gvec{\lambda}-\gvec{\lambda}^{(t)}) \\+\vec{g}^T(\gvec{\lambda}^{(t)})(\gvec{\lambda}-\gvec{\lambda}^{(t)}) + d(\gvec{\lambda}^{(t)}).
\end{multline}
This follows the idea of Newton methods, where the gradient and Hessian of the function are used within the approximation. However, due to the nonsmoothness of the dual function, the gradient and Hessian are not defined for each value of the dual variable. Instead, the gradient is replaced in eq. (\ref{QNDA_Dual}) by the subgradient and the Hessian is approximated by the matrix $ \vec{B}^{(t)} $. The approximated Hessian can be updated in each iteration using a Broyden-Fletcher-Goldfarb-Shanno (BFGS) update,
\begin{equation}\label{BFGSUpdate}
	\vec{B}^{(t)} = \vec{B}^{(t-1)} + \frac{\vec{y}^{(t)}\vec{y}^{(t),T}}{\vec{y}^{(t),T}\vec{s}^{(t)}} - \frac{\vec{B}^{(t-1)}\vec{s}^{(t)}\vec{s}^{(t),T}\vec{B}^{(t-1),T}}{\vec{s}^{(t),T}\vec{B}^{(t-1)}\vec{s}^{(t)}},
\end{equation}
where
\begin{equation}\label{Stepsize}
	\vec{s}^{(t)} \coloneqq \gvec{\lambda}^{(t)} - \gvec{\lambda}^{(t-1)}
\end{equation}
is the variation of the dual variables and 
\begin{equation}\label{SubgradientDifference}
	\vec{y}^{(t)} = \vec{g}(\gvec{\lambda}^{(t)}) - \vec{g}(\gvec{\lambda}^{(t-1)})
\end{equation}
is the variation of the subgradients.\\
The approximated dual function $ d_B(\gvec{\lambda}) $ is differentiable, while the actual dual function is nonsmooth. This can lead to significant approximation errors and poor update directions. This issue can be addressed by utilizing the same information as in the BTM algorithm. However, instead of using the bundle to construct an over-approximator of the dual function, it is used to further constrain the update of the dual variables,
\begin{multline}\label{Bundle_Cut}
	d_B^{(t)}(\gvec{\lambda}^{(t+1)}) \le d(\gvec{\lambda}^{(l)}) + \vec{g}^T(\gvec{\lambda}^{(l)}) (\gvec{\lambda}^{(t+1)} - \gvec{\lambda}^{(l)}),\\\forall l \in \{t-\tau+1,\dots,t\}.
\end{multline}

Constraints (\ref{Bundle_Cut}) are derived from the definition of the subgradient (\ref{Subgradient_Condition}). A violation of these constraints would indicate that the updated dual variables $ \gvec{\lambda}^{(t+1)} $ are outside the range of validity of the approximated dual function. These constraints are referred to as bundle cuts and they can be summarized as
\begin{multline}\label{BC_Summary}
		\mathcal{BC}^{(t)} = \{\gvec{\lambda}\in\vbody{R}{\vec{b}}\vert\;d_B^{(t)}(\gvec{\lambda}) \le d(\gvec{\lambda}^{(l)}) + \vec{g}^T(\gvec{\lambda}^{(l)}) (\gvec{\lambda} - \gvec{\lambda}^{(l)}),\\	\forall l \in \{t-\tau+1,\dots,t\}\}.
\end{multline}
In the QNDA algorithm, the dual variables are updated in each iteration by solving the optimization problem
	\begin{subequations}\label{QNDA_Update}
		\begin{align}
			\gvec{\lambda}^{(t+1)} = \text{arg}\underset{\gvec{\lambda}}{\max}\;&d_B^{(t)}(\gvec{\lambda}),\\
			\text{s.\,t. }&\Vert \gvec{\lambda} - \gvec{\lambda}^{(t)}\Vert_2^2 \le \alpha^{(t)},\\
			&\gvec{\lambda} \in \mathcal{BC}^{(t)}.
		\end{align}
	\end{subequations}
To avoid too aggressive update steps the same trust region (\ref{QNDA_Update}b) as in the BTM algorithm is used.
\subsection{Primal heuristics}
The following sections provide some additional heuristics related to the primal optimization problem (\ref{Clustering_MILP_Consenus}), namely an averaging heuristic used to obtain feasible primal solutions, and the addition of symmetry-breaking constraints to the clustering problem.
\subsubsection{Averaging heuristic}\label{Sec:Avg_Heuristic}
\begin{figure*}[h!]
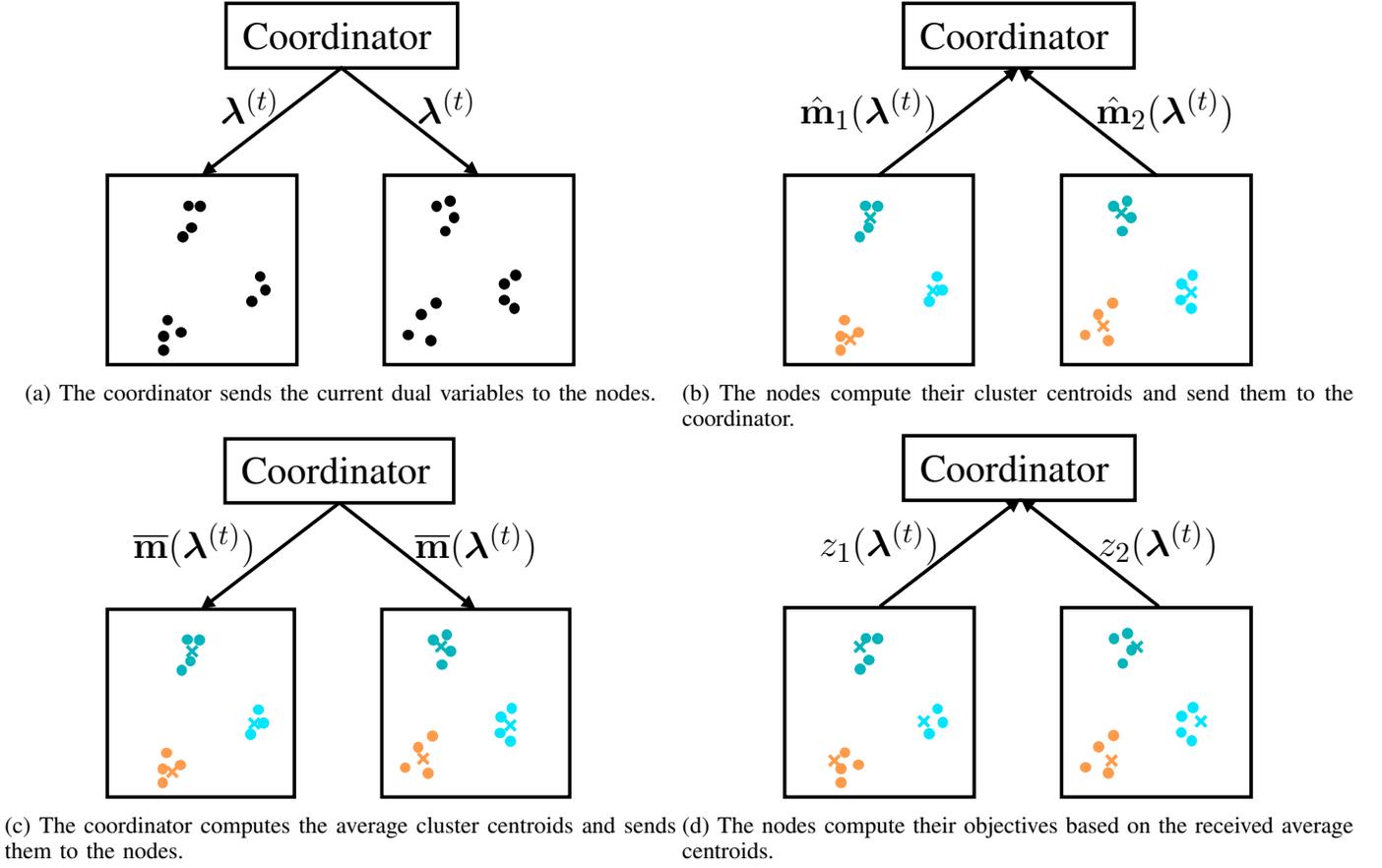

	\begin{subfigure}[t]{0.5\linewidth}
		\centering
		\fontsize{15pt}{3pt}\selectfont
		\includesvg[width = 0.7\linewidth]{Figures/Clustering/Distr_Clustering_1}
		\caption{The coordinator sends the current dual variables to the nodes.}
		\label{fig:distr_clustering_1}
	\end{subfigure}
	\begin{subfigure}[t]{0.5\linewidth}
		\centering
		\fontsize{15pt}{3pt}\selectfont
		\includesvg[width = 0.7\linewidth]{Figures/Clustering/Distr_Clustering_2}
		\caption{The nodes compute their cluster centroids and send them to the coordinator.}
		\label{fig:distr_clustering_2}
	\end{subfigure}
	
	\begin{subfigure}[t]{0.5\linewidth}
		\centering
		\fontsize{15pt}{3pt}\selectfont
		\includesvg[width = 0.7\linewidth]{Figures/Clustering/Distr_Clustering_3}
		\caption{The coordinator computes the average cluster centroids and sends them to the nodes.}
		\label{fig:distr_clustering_3}
	\end{subfigure}
	\begin{subfigure}[t]{0.5\linewidth}
		\centering
		\fontsize{15pt}{3pt}\selectfont
		\includesvg[width = 0.7\linewidth]{Figures/Clustering/Distr_Clustering_4}
		\caption{The nodes compute their objectives based on the received average centroids.}
		\label{fig:distr_clustering_4}
	\end{subfigure}
	\caption{Communication process between the coordinator and the nodes in iteration $ t $.}
	\label{fig:distr_clustering}
\end{figure*}

The $ K $-means clustering problem involves integrality constraints and is therefore nonconvex. While the (optimal) value of the dual function \eqref{Clustering_Dual} provides a lower bound on the optimal value of the primal problem \eqref{Clustering_MILP_Consenus}, the feasibility of the primal problem is not guaranteed upon the convergence of a dual decomposition-based algorithm, i.e., the consensus constraints may not be satisfied. Nevertheless, in the case of $ K $-means clustering it is straightforward to compute a feasible primal solution using an averaging step. In each iteration $ t $ of a dual decomposition-based algorithm the coordinator communicates the dual variables $ \gvec{\lambda}^{(t)} $ to the nodes. The nodes in turn solve their clustering problems and communicate their computed cluster centroids $ \hat{\vec{m}}_i(\gvec{\lambda}^{(t)}) $ to the coordinator. Based on this response the coordinator can compute the average of the primal variables, i.e., the average cluster centroids,
\begin{equation}\label{Avg_Centroids}
	\overline{\vec{m}}_k(\gvec{\lambda}^{(t)}) = \frac{1}{N_s} \sum_{i\in\mathcal{I}} \vec{m}_{ik}(\gvec{\lambda}^{(t)})
\end{equation}
which are then communicated back to the nodes. Using the mean cluster centroids the nodes can compute their resulting primal objective value
\begin{equation}\label{Node_Objective}
	z_i(\gvec{\lambda}^{(t)}) = \sum_{j\in\mathcal{J}} \underset{k\in\mathcal{K}}{\min}\Vert \vec{y}_j - \overline{\vec{m}}_k(\gvec{\lambda}^{(t)}) \Vert_2^2.
\end{equation}
The primal objective value can be used to compute the relative duality gap in each iteration,
\begin{equation}\label{Clustering_DG}
	\text{rel. DG} = 100\cdot\left(1-\frac{d(\gvec{\lambda}^{(t)})}{\sum_{i\in\mathcal{I}}z_i(\gvec{\lambda}^{(t)}) }\right).
\end{equation}
Since the value of the dual function provides a lower bound on the optimal primal objective value the relative duality gap can be used to assess the distance of a found solution to the global optimum. The entire communication process between the coordinator and the nodes is illustrated in Fig.~\ref{fig:distr_clustering}. Note that the average cluster centroids are only used to compute the duality gap. They do not influence the update of the dual variables.
\subsubsection{Symmetry breaking constraints}
\begin{figure}
	\centering
	\fontsize{8pt}{3pt}\selectfont
	\includesvg[width = 0.85\linewidth]{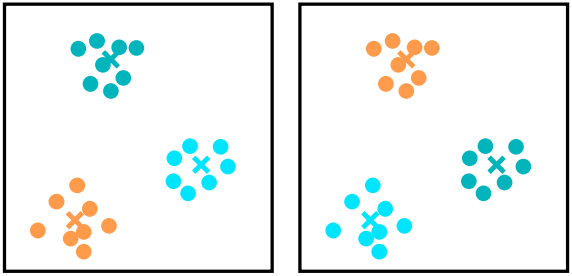}
	\caption[Example of symmetric clustering solutions.]{Example of symmetric clustering solutions. In the two cases, the data points are assigned to different clusters without affecting the objective function.}
	\label{fig:symmetric_clustering}
\end{figure}

The clustering problem \eqref{Clustering_MILP_Central} is highly symmetric, i.e., it contains solutions with the same objective values. This is because the index assigned to a cluster does not influence the objective function. Fig.~\ref{fig:symmetric_clustering} illustrates the situation of two symmetric solutions.
This symmetry can lead to problems for the averaging heuristic presented in the previous section, as the computed cluster centroids of a single node can switch from one iteration to the next. For instance, while some points are assigned to cluster $ k $ in iteration $ t $, they could be assigned to cluster $ k' $ in iteration $ t+1 $ by switching the centroids of clusters $ k $ and $ k' $ without affecting the objective.\\
To prevent this behavior symmetry breaking constraints are added to the optimization problems of the nodes. In the first iteration, one of the nodes acts as the reference node, providing reference centroids $ \overline{\vec{m}}_k^{\text{ref}} $. In the subsequent iterations the quadratic constraint
\begin{equation}\label{Symmetry_Constraint}
	\Vert \vec{m}_{ik} - \overline{\vec{m}}_k^{\text{ref}} \Vert_2^2 \le \Vert \vec{m}_{ik'} - \overline{\vec{m}}_k^{\text{ref}} \Vert_2^2, \forall k,k'\in\mathcal{K},
\end{equation}
is added to each node $ i $. This ensures that cluster $ k $ of each node $ i $ will be the one closest to the reference centroid $ \overline{\vec{m}}_k^{\text{ref}} $. The choice of the node which provides the reference centroid can be performed arbitrarily, as it does not affect the optimization of the other nodes. Furthermore, the added constraint also does not affect the optimal objective value while rendering all symmetric solutions, except for one, infeasible.
\section{Numerical analysis of distributed clustering problems}
\begin{table}[h!]
	\centering
	\caption{Parameter settings of the distributed optimization algorithms for the clustering benchmark problems.}
	\begin{tabular}{l l l l}
		\hline
		& \textbf{Value} &  \textbf{Description} & \textbf{Algorithms} \\
		\hline
		$ \gvec{\lambda}^{(0)} $ &  $\vec{0}$ & initial dual variables & All\\
		$ \alpha^{(0)} $&  $0.5$ & initial step size/trust & All\\
  		  &   & region parameter & \\ 
		
		$ t_{\max} $& $150$ &   maximum number of & All\\
		&   & iterations & \\
		$ \epsilon_p $& $10^{-2}$ &  primal residual convergence  & All\\
		&   & tolerance& \\
		$ \epsilon_{DG} $ & $ 0.25\; \% $ & relative duality gap  & All\\
		&   & tolerance& \\
		$ \tau $& $50 $ &  allowed age of data  & BTM, QNDA\\
		&   & points& \\
		$ \vec{B}^{(0)} $ & $ -\vec{I} $ & initial approximated  & QNDA\\
				
		& & Hessian & \\
        & &  & \\
		\hline
	\end{tabular}
	\label{tab:Parameter_Table_Clustering}
\end{table}

The dual decomposition-based distributed clustering approach was evaluated on a set of benchmark problems of varying sizes. The data for each benchmark problem was generated randomly. First, initial cluster centroids $ \vec{m}_{k}^0 $ were generated, with $  [\vec{m}_{k}^0]_l \in \mathcal{U}_c(-1,1),\;l = 1,\dots,n_{\vec{y}} $.  Then, for each cluster $ k $ five random data points were added within a radius of 0.5 from the generated centroid. The parameters of the benchmark problems were varied as follows:
\pagebreak
\begin{align*}
	&\text{Number of nodes: } N_s \in \{2,3,4\},\\
	&\text{Number of dimensions: }n_{\vec{y}} \in \{2,3,4\},\\
	&\text{Number of clusters: }K \in \{3,4\}.
\end{align*}
Five benchmark problems were generated for each combination of nodes, dimensions, and clusters, resulting in a total of 90 benchmark problems. A benchmark problem is characterized by its number of nodes, dimension of the data, and number of clusters. For instance, problem  $ \text{3N2D4K}_{\text{5}} $ is the 5th benchmark problem comprised of 3 nodes with 2-dimensional data sorted into 4 clusters.\\
The benchmark problems were solved using the subgradient method, the bundle trust method, and the quasi-Newton dual ascent algorithm.

The initial step size (subgradient method)/ trust region (BTM, QNDA) parameter was set to $ \alpha^{(0)} = 0.5 $ and varied according to
\begin{equation}\label{Step_Size_Clustering}
	\alpha^{(t)} = \alpha^{(0)}/\sqrt{t}.
\end{equation}
The size of the bundle for BTM and QNDA was set to $ \tau = 50 $ points. All algorithms were initialized with $ \gvec{\lambda}^{(0)} = \vec{0} $ and the initial approximated Hessian of the QNDA algorithm was set to the negative identity matrix. The algorithms were terminated either when the Euclidean norm of the primal residual
\begin{equation}\label{key}
	\Vert \vec{w}_p\Vert_2 = \left\Vert \sum_{i\in\mathcal{I}}\vec{A}_i\hat{\vec{m}}_i \right\Vert_2,
\end{equation}
i.e., the violation of the consensus constraints, lied below a threshold of $ \epsilon_p = 10^{-2} $ or when the relative duality gap \eqref{Clustering_DG} reached a value of $ \epsilon_{DG} = 0.25\;\% $. The used parameters for the different algorithms are summarized in Tab.~\ref{tab:Parameter_Table_Clustering}. The MIQCP clustering problems of all nodes were solved using the commercial solver Gurobi \cite{gurobi} and the total computation time was computed as
\begin{equation}\label{Comp_Time}
	T_{\text{comp}} = N_{\text{iter}}\cdot T_{\text{comm}} + \sum_{t = 1}^{N_{\text{iter}}}  (T_{\text{update}}^{(t)} +    \underset{i\in\mathcal{I}}{\max}\;T_{\text{sub},i}^{(t)}),
\end{equation}
where $ N_{\text{iter}} $ is the number of required iterations, $ T_{\text{comm}}=800\;ms $ is the required communication time between the coordinator and the subproblems, which is assumed to be constant, $ T_{\text{update}}^{(t)} $ is the time required by the coordinator to update the dual variables in iteration $ t $ and $ T_{\text{sub},i}^{(t)} $ is the solution time for the clustering problem of node $ i $ in itereation $ t $.

\begin{table}[h!]
	\centering
	\caption[Summary of the results for the distributed optimization of the clustering benchmark problems.]{Summary of the results for the distributed optimization of the clustering benchmark problems, $ \overline{t} $: mean number of iterations until termination, $ \overline{\text{rel. DG}}$: mean relative duality gap upon termination (in \%), $ \overline{T_\text{comp}} $: mean computation time (in s).}
\begin{tabular}{l|c c c }
	Algorithm & $ \overline{t} $ & $ \overline{\text{rel. DG}} $ & $\overline{T_\text{comp}}$  \\
	\hline
	SG &${136.75}$	&${2.27}$	&	${996.28}$\\
	
	BTM &${57.44}$	&${1.86}$	&${515.77}$  \\
	
	QNDA &$\vec{54.48}$	&$\mathbf{1.81}$	&$\vec{483.22}$ \\
	
\end{tabular}
\label{tab:Results_Clustering_Summary}
\end{table}

The results for the clustering benchmarks are summarized in Tab.~\ref{tab:Results_Clustering_Summary}. Out of the examined algorithms, QNDA shows the best performance in terms of the required number of iterations and computation time as well as in terms of the achieved relative duality gap. The BTM algorithm shows similar performance in terms of the number of iterations and the achieved duality gap. However, in the case of distributed clustering, each iteration is costly due to the underlying MIQCP problems. Therefore, a slight performance increase in the number of iterations results in a substantial performance increase in terms of computation times. More detailed results for the clustering benchmarks are summarized in Tab.~\ref{tab:Results_Clustering} in the appendix.
\begin{figure}[h!]
	\centering
	\fontsize{7pt}{3pt}\selectfont
	\includesvg[width = 0.9\linewidth]{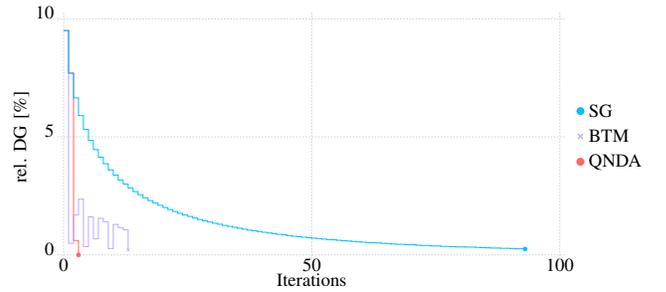}
	\caption[Evolution of the relative duality gap for benchmark problem $ \text{2N2D4K}_{\text{4}} $.]{Evolution of the relative duality gap for benchmark problem $ \text{2N2D4K}_{\text{3}}$.}
	\label{fig:rel_DG_2N2D4K_3}	
\end{figure}

\begin{figure}[h]
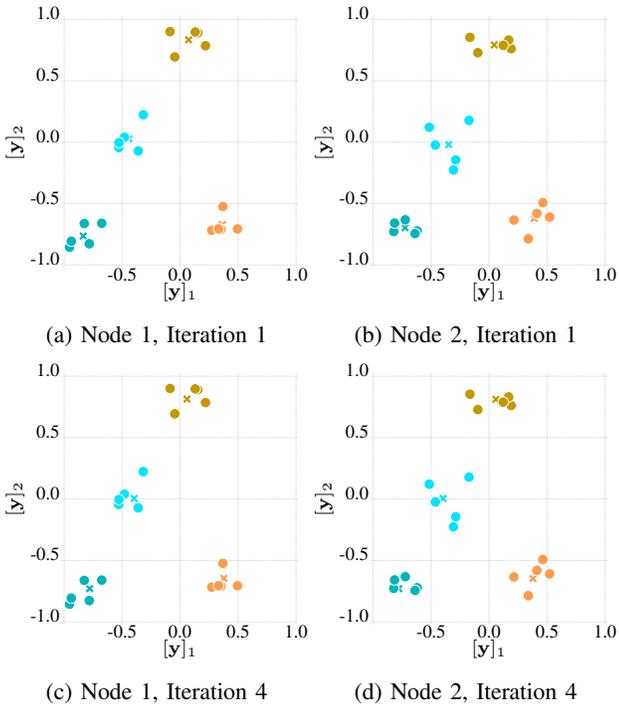

	\begin{subfigure}[b]{0.45\linewidth}
		\centering
		\fontsize{7pt}{3pt}\selectfont
		\includesvg[width = \linewidth]{Figures/Clustering/2N2D4K_3/Node_1Iter_1}
		\caption{Node 1, Iteration 1}
		\label{fig:Node_1_It_1}
	\end{subfigure}
	\begin{subfigure}[b]{0.45\linewidth}
		\centering
		\fontsize{7pt}{3pt}\selectfont
		\includesvg[width = \linewidth]{Figures/Clustering/2N2D4K_3/Node_2Iter_1}
		\caption{Node 2, Iteration 1}
		\label{fig:Node_2_It_1}
	\end{subfigure}
	
	\vspace{2mm}
	\begin{subfigure}[b]{0.45\linewidth}
		\centering
		\fontsize{7pt}{3pt}\selectfont
		\includesvg[width = \linewidth]{Figures/Clustering/2N2D4K_3/Node_1Iter_4}
		\caption{Node 1, Iteration 4}
		\label{fig:Node_1_It_4}
	\end{subfigure}
	\begin{subfigure}[b]{0.45\linewidth}
		\centering
		\fontsize{7pt}{3pt}\selectfont
		\includesvg[width = \linewidth]{Figures/Clustering/2N2D4K_3/Node_2Iter_4}
		\caption{Node 2, Iteration 4}
		\label{fig:Node_2_It_4}
	\end{subfigure}
	\caption{Exemplary clusters in different iterations of the QNDA algorithm for benchmark problem $ \text{2N2D4K}_{\text{3}} $.}
	\label{fig:2N2D4K_3_Iterations}
\end{figure}

Fig.~\ref{fig:rel_DG_2N2D4K_3} shows the evolution of the relative duality gap for benchmark problem $ \text{2N2D4K}_{\text{3}}$. The subgradient method converges rather slowly. In comparison, the BTM and QNDA algorithms exhibit a faster rate of convergence. Between these two algorithms, BTM exhibits an oscillatory behavior before converging. In contrast, the QNDA algorithm does not exhibit oscillations and therefore converges earlier. Additionally, it should be noted that the QNDA algorithm achieves a relative duality gap of $ 0\; \%$, i.e., it converges to a proven global optimum.

{Fig.~\ref{fig:2N2D4K_3_Iterations} further illustrates the results. Fig.~\ref{fig:Node_1_It_1} and \ref{fig:Node_2_It_1} show the results of the clustering in the first iteration, i.e., the individual global optima. Fig.~\ref{fig:Node_1_It_4} and \ref{fig:Node_2_It_4} depict the solutions upon convergence of the QNDA algorithm. Each node computes the same cluster centroids corresponding to the globally optimal solution with respect to the entire data set, but not to the individual data sets. It is therefore possible to compute a global model locally in each node while only accessing local data.}

\section{Comparison to the central solution}
\begin{figure}[h]
	\centering
	\fontsize{7pt}{3pt}\selectfont
	\includesvg[width = 0.9\linewidth]{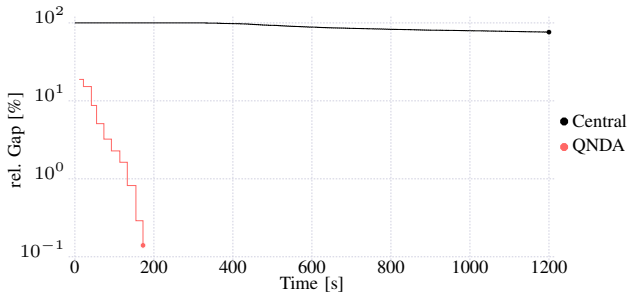}
	\caption{Evolution of the relative duality gap of QNDA compared to the relative integrality gap of the central solution using Gurobi for benchmark problem $ \text{4N4D4K}_{3} $.}
	\label{fig:central_clustering}	
\end{figure}
As shown in the previous section, solving the MIQCP clustering problems is computationally expensive. This is due to the weak integer relaxation of problem \eqref{Clustering_MILP}, which means that the solution of the relaxed problem within the branch-and-bound algorithm is far away from the integer solution. This results in slow-moving relative integrality gaps and slow convergence of the solution algorithm. While the main motivation of the distributed clustering approach is the training of a global model without the exchange of local data, it can also be used to efficiently solve larger clustering problems. Fig.~\ref{fig:central_clustering}	 depicts the evolution of the relative duality gap of the QNDA algorithm as well as the evolution of the relative integrality gap of Gurobi for the complete data set of benchmark problem $ \text{4N4D4K}_{3} $. The clustering problems of the individual nodes were solved sequentially in the case of QNDA, also using Gurobi. While the relative gap of the central solution improves very slowly, the QNDA algorithm quickly converges to a solution close to the global optimum. Note, that both relative gaps prove a worst-case distance to the global optimum. Hence, decomposing a large clustering problem into smaller subproblems and coordinating the solutions via a distributed optimization algorithm can offer significant performance improvements compared to a central solution.
\section{Conclusion}
This paper demonstrated how dual decomposition-based distributed optimization can be applied to the solution of clustering problems. The approach ensures privacy, i.e., enables federated learning, as each node only has access to its local data. A global model can still be obtained by coordinating the solutions of the individual clustering problems. Numerical tests on a large set of benchmark problems demonstrated that the QNDA algorithm outperforms the subgradient method and the BTM algorithm. Furthermore, the distributed optimization approach exhibited superior performance compared to a central solution approach. In the future, the developed algorithms can also be applied to other federated learning problems, like the distributed training of support vector machines.
\bibliography{main}

\begin{thebibliography}{10}
\providecommand{\url}[1]{#1}
\csname url@samestyle\endcsname
\providecommand{\newblock}{\relax}
\providecommand{\bibinfo}[2]{#2}
\providecommand{\BIBentrySTDinterwordspacing}{\spaceskip=0pt\relax}
\providecommand{\BIBentryALTinterwordstretchfactor}{4}
\providecommand{\BIBentryALTinterwordspacing}{\spaceskip=\fontdimen2\font plus
\BIBentryALTinterwordstretchfactor\fontdimen3\font minus
  \fontdimen4\font\relax}
\providecommand{\BIBforeignlanguage}[2]{{%
\expandafter\ifx\csname l@#1\endcsname\relax
\typeout{** WARNING: IEEEtran.bst: No hyphenation pattern has been}%
\typeout{** loaded for the language `#1'. Using the pattern for}%
\typeout{** the default language instead.}%
\else
\language=\csname l@#1\endcsname
\fi
#2}}
\providecommand{\BIBdecl}{\relax}
\BIBdecl

\bibitem{Peteiro2013}
D.~{Peteiro-Barral} and B.~{Guijarro-Berdi{\~n}as}, ``A survey of methods for
  distributed machine learning,'' \emph{Progress in Artificial Intelligence},
  vol.~2, pp. 1--11, 2013.

\bibitem{Forero2010}
P.~Forero, A.~Cano, and G.~Giannakis, ``Consensus-based distributed support
  vector machines.'' \emph{Journal of Machine Learning Research}, vol.~11,
  no.~5, 2010.

\bibitem{Forero2011}
------, ``Distributed clustering using wireless sensor networks,'' \emph{IEEE
  Journal of Selected Topics in Signal Processing}, vol.~5, no.~4, pp.
  707--724, 2011.

\bibitem{Georgopoulos2014}
L.~Georgopoulos and M.~Hasler, ``Distributed machine learning in networks by
  consensus,'' \emph{Neurocomputing}, vol. 124, pp. 2--12, 2014.

\bibitem{Tsianos2012}
K.~Tsianos, S.~Lawlor, and M.~Rabbat, ``Consensus-based distributed
  optimization: Practical issues and applications in large-scale machine
  learning,'' in \emph{50th Annual Allerton Conference on Communication,
  Control, and Computing}.\hskip 1em plus 0.5em minus 0.4em\relax IEEE, 2012,
  pp. 1543--1550.

\bibitem{Nedic2020}
A.~Nedi{\'c}, ``Distributed gradient methods for convex machine learning
  problems in networks: Distributed optimization,'' \emph{IEEE Signal
  Processing Magazine}, vol.~37, no.~3, pp. 92--101, 2020.

\bibitem{Verbraeken2020}
J.~Verbraeken, M.~Wolting, J.~Katzy, J.~Kloppenburg, T.~Verbelen, and
  J.~Rellermeyer, ``A survey on distributed machine learning,'' \emph{ACM
  Computing Surveys}, vol.~53, no.~2, pp. 1--33, 2020.

\bibitem{Konevcny2016}
J.~Kone{\v{c}}n{\`y}, B.~McMahan, D.~Ramage, and P.~Richt{\'a}rik, ``Federated
  optimization: Distributed machine learning for on-device intelligence,''
  \emph{arXiv preprint arXiv:1610.02527}, 2016.

\bibitem{Ghosh2020}
A.~Ghosh, J.~Chung, D.~Yin, and K.~Ramchandran, ``An efficient framework for
  clustered federated learning,'' \emph{Advances in Neural Information
  Processing Systems}, vol.~33, pp. 19\,586--19\,597, 2020.

\bibitem{Mcmahan2017}
B.~McMahan, E.~Moore, D.~Ramage, S.~Hampson, and B.~A. y~Arcas,
  ``Communication-efficient learning of deep networks from decentralized
  data,'' in \emph{Artificial Intelligence and Statistics}.\hskip 1em plus
  0.5em minus 0.4em\relax PMLR, 2017, pp. 1273--1282.

\bibitem{Hegiste2022}
V.~Hegiste, T.~Legler, and M.~Ruskowski, ``Application of federated machine
  learning in manufacturing,'' in \emph{2022 International Conference on
  Industry 4.0 Technology (I4Tech)}.\hskip 1em plus 0.5em minus 0.4em\relax
  IEEE, 2022, pp. 1--8.

\bibitem{Antunes2022}
R.~Antunes, C.~Andr{\'e}\space{}da\space{}Costa, A.~K{\"u}derle, I.~Yari, and
  B.~Eskofier, ``Federated learning for healthcare: Systematic review and
  architecture proposal,'' \emph{ACM Transactions on Intelligent Systems and
  Technology (TIST)}, vol.~13, no.~4, pp. 1--23, 2022.

\bibitem{Lim2020}
W.~Lim, N.~Luong, D.~Hoang, Y.~Jiao, Y.-C. Liang, Q.~Yang, D.~Niyato, and
  C.~Miao, ``Federated learning in mobile edge networks: A comprehensive
  survey,'' \emph{IEEE Communications Surveys \& Tutorials}, vol.~22, no.~3,
  pp. 2031--2063, 2020.

\bibitem{Jiang2020}
J.~Jiang, B.~Kantarci, S.~Oktug, and T.~Soyata, ``Federated learning in smart
  city sensing: Challenges and opportunities,'' \emph{Sensors}, vol.~20,
  no.~21, p. 6230, 2020.

\bibitem{Li2020}
L.~Li, Y.~Fan, M.~Tse, and K.-Y. Lin, ``A review of applications in federated
  learning,'' \emph{Computers \& Industrial Engineering}, vol. 149, p. 106854,
  2020.

\bibitem{Liu2022}
J.~Liu, J.~Huang, Y.~Zhou, X.~Li, S.~Ji, H.~Xiong, and D.~Dou, ``From
  distributed machine learning to federated learning: A survey,''
  \emph{Knowledge and Information Systems}, vol.~64, no.~4, pp. 885--917, 2022.

\bibitem{Chamikara2021}
M.~Chamikara, P.~Bertok, I.~Khalil, D.~Liu, and S.~Camtepe, ``Privacy
  preserving distributed machine learning with federated learning,''
  \emph{Computer Communications}, vol. 171, pp. 112--125, 2021.

\bibitem{Kulkarni2020}
V.~Kulkarni, M.~Kulkarni, and A.~Pant, ``Survey of personalization techniques
  for federated learning,'' in \emph{4th World Conference on Smart Trends in
  Systems, Security and Sustainability (WorldS4)}.\hskip 1em plus 0.5em minus
  0.4em\relax IEEE, 2020, pp. 794--797.

\bibitem{Tan2022}
A.~Tan, H.~Yu, L.~Cui, and Q.~Yang, ``Towards personalized federated
  learning,'' \emph{IEEE Transactions on Neural Networks and Learning Systems},
  2022.

\bibitem{Gambella2021}
C.~Gambella, B.~Ghaddar, and J.~{Naoum-Sawaya}, ``Optimization problems for
  machine learning: A survey,'' \emph{European Journal of Operational
  Research}, vol. 290, no.~3, pp. 807--828, 2021.

\bibitem{Dhanachandra2015}
N.~Dhanachandra, K.~Manglem, and Y.~Chanu, ``Image segmentation using k-means
  clustering algorithm and subtractive clustering algorithm,'' \emph{Procedia
  Computer Science}, vol.~54, pp. 764--771, 2015.

\bibitem{Kansal2018}
T.~Kansal, S.~Bahuguna, V.~Singh, and T.~Choudhury, ``Customer segmentation
  using k-means clustering,'' in \emph{International Conference on
  Computational Techniques, Electronics and Mechanical Systems (CTEMS)}.\hskip
  1em plus 0.5em minus 0.4em\relax IEEE, 2018, pp. 135--139.

\bibitem{Rahimi2019}
K.~{Rahimi-Adli}, P.~Schiermoch, B.~Beisheim, S.~Wenzel, and S.~Engell, ``A
  model identification approach for the evaluation of plant efficiency,'' in
  \emph{Computer Aided Chemical Engineering}.\hskip 1em plus 0.5em minus
  0.4em\relax Elsevier, 2019, vol.~46, pp. 913--918.

\bibitem{Dennis2021heterogeneity}
D.~K. Dennis, T.~Li, and V.~Smith, ``Heterogeneity for the win: One-shot
  federated clustering,'' in \emph{International Conference on Machine
  Learning}.\hskip 1em plus 0.5em minus 0.4em\relax PMLR, 2021, pp. 2611--2620.

\bibitem{Lloyd1982}
S.~Lloyd, ``Least squares quantization in pcm,'' \emph{IEEE Transactions on
  Information Theory}, vol.~28, no.~2, pp. 129--137, 1982.

\bibitem{Kumar2020}
H.~Kumar, V.~Karthik, and M.~Nair, ``Federated k-means clustering: A novel edge
  ai based approach for privacy preservation,'' in \emph{2020 IEEE
  International Conference on Cloud Computing in Emerging Markets
  (CCEM)}.\hskip 1em plus 0.5em minus 0.4em\relax IEEE, 2020, pp. 52--56.

\bibitem{Li2022}
\BIBentryALTinterwordspacing
S.~Li, S.~Hou, B.~Buyukates, and S.~Avestimehr, ``Secure federated
  clustering,'' \emph{arXiv}, 2022. [Online]. Available:
  \url{https://arxiv.org/abs/2205.15564}
\BIBentrySTDinterwordspacing

\bibitem{Stallmann2022}
\BIBentryALTinterwordspacing
M.~Stallmann and A.~Wilbik, ``Towards federated clustering: A federated fuzzy $
  c $-means algorithm ({FFCM}),'' \emph{arXiv}, 2022. [Online]. Available:
  \url{https://arxiv.org/abs/2201.07316}
\BIBentrySTDinterwordspacing

\bibitem{Pedrycz2021}
W.~Pedrycz, ``Federated {FCM}: clustering under privacy requirements,''
  \emph{IEEE Transactions on Fuzzy Systems}, vol.~30, no.~8, pp. 3384--3388,
  2021.

\bibitem{Wang2022federated}
Y.~Wang, M.~Jia, N.~Gao, L.~Von~Krannichfeldt, M.~Sun, and G.~Hug, ``Federated
  clustering for electricity consumption pattern extraction,'' \emph{IEEE
  Transactions on Smart Grid}, vol.~13, no.~3, pp. 2425--2439, 2022.

\bibitem{Achterberg2013}
T.~Achterberg and R.~Wunderling, ``Mixed integer programming: Analyzing 12
  years of progress,'' \emph{Facets of Combinatorial Optimization: Festschrift
  for Martin Gr{\"o}tschel}, pp. 449--481, 2013.

\bibitem{Koch2022}
T.~Koch, T.~Berthold, J.~Pedersen, and C.~Vanaret, ``Progress in mathematical
  programming solvers from 2001 to 2020,'' \emph{EURO Journal on Computational
  Optimization}, vol.~10, p. 100031, 2022.

\bibitem{Aloise2012}
D.~Aloise, P.~Hansen, and L.~Liberti, ``An improved column generation algorithm
  for minimum sum-of-squares clustering,'' \emph{Mathematical Programming},
  vol. 131, pp. 195--220, 2012.

\bibitem{Bagirov2006}
A.~Bagirov and J.~Yearwood, ``A new nonsmooth optimization algorithm for
  minimum sum-of-squares clustering problems,'' \emph{European Journal of
  Operational Research}, vol. 170, no.~2, pp. 578--596, 2006.

\bibitem{Karmitsa2017}
N.~Karmitsa, A.~Bagirov, and S.~Taheri, ``New diagonal bundle method for
  clustering problems in large data sets,'' \emph{European Journal of
  Operational Research}, vol. 263, no.~2, pp. 367--379, 2017.

\bibitem{Everett.1963}
H.~Everett, ``Generalized {L}agrange multiplier method for solving problems of
  optimum allocation of resources,'' \emph{Operations Research}, no. 11 (3),
  pp. 399--417, 1963.

\bibitem{Nocedal2006}
J.~Nocedal and S.~Wright, \emph{Numerical optimization}.\hskip 1em plus 0.5em
  minus 0.4em\relax Springer Science \& Business Media, 2006.

\bibitem{Yfantis2023EURO}
V.~Yfantis, S.~Wenzel, A.~Wagner, M.~Ruskowski, and S.~Engell, ``Hierarchical
  distributed optimization of constraint-coupled convex and mixed-integer
  programs using approximations of the dual function,'' \emph{EURO Journal on
  Computational Optimization}, vol.~11, p. 100058, 2023.

\bibitem{Shor2012minimization}
N.~Shor, \emph{Minimization methods for non-differentiable functions}.\hskip
  1em plus 0.5em minus 0.4em\relax Springer Science \& Business Media, 2012,
  vol.~3.

\bibitem{Bertsekas.1999}
D.~P. Bertsekas, \emph{Nonlinear programming}.\hskip 1em plus 0.5em minus
  0.4em\relax {Athena Scientific}, 1999.

\bibitem{Makela.2002}
M.~M{\"a}kel{\"a}, ``Survey of bundle methods for nonsmooth optimization,''
  \emph{Optimization Methods and Software}, vol.~17, no.~1, pp. 1--29, 2002.

\bibitem{Bagirov2014}
A.~Bagirov, N.~Karmitsa, and M.~M{\"a}kel{\"a}, \emph{Introduction to Nonsmooth
  Optimization: {T}heory, {P}ractice and {S}oftware}.\hskip 1em plus 0.5em
  minus 0.4em\relax Springer, 2014.

\bibitem{Le2007}
Q.~Le, A.~Smola, and S.~Vishwanathan, ``Bundle methods for machine learning,''
  \emph{Advances in Neural Information Processing Systems}, vol.~20, 2007.

\bibitem{Yfantis2022}
V.~Yfantis and M.~Ruskowski, ``A hierarchical dual decomposition-based
  distributed optimization algorithm combining quasi-{N}ewton steps and bundle
  methods,'' in \emph{30th Mediterranean Conference on Control and Automation
  (MED)}.\hskip 1em plus 0.5em minus 0.4em\relax IEEE, 2022, pp. 31--36.

\bibitem{gurobi}
\BIBentryALTinterwordspacing
{Gurobi Optimization, LLC}, ``{Gurobi Optimizer Reference Manual},'' 2023.
  [Online]. Available: \url{https://www.gurobi.com}
\BIBentrySTDinterwordspacing

\end{thebibliography}
\bibliographystyle{IEEEtran}
\appendix
\subsection{Benchmark Problems}
All benchmark problems can be found under:
\url{https://github.com/VaYf/Clustering-Benchmark-Problems}
%\subsection{Results for the clustering benchmark problems}
%The results for the clustering benchmarks are summarized in Tab. \ref{tab:Results_Clustering}.
\onecolumn
\begin{table*}[ht]
 \centering
		\caption[Results for the distributed optimization of the clustering benchmark problems.]{Results for the distributed optimization of the clustering benchmark problems, $ \overline{t} $: mean number of performed iterations, $ \overline{\text{rel. DG}} $: mean relative duality gap (in \%), $\overline{T_\text{comp}}$: mean computation time (in s).}\label{tab:Results_Clustering}
		%\endfirsthead
		%\hline \multicolumn{7}{r}{{Continued on next page}} \\ \hline
		%\endfoot
		%\endlastfoot
    \begin{tabular}{c | c c c | c c c}
		\multicolumn{1}{c}{}&\multicolumn{3}{|c}{SG} & \multicolumn{3}{|c}{BTM} \\  
		Clustering & $ \overline{t} $ &   $ \overline{\text{rel. DG}} $  & $\overline{T_\text{comp}}$ & $ \overline{t} $ &   $ \overline{\text{rel. DG}} $  & $\overline{T_\text{comp}}$ \\
		\hline
		\textbf{Mean}	& $\vec{136.75}$	& $\vec{2.27}$	& $\vec{996.28}$	& $\vec{57.44}$	& $\vec{1.86}$	& $\vec{515.77}$\\
		\hline
		2N2D3K	& $126.0$	& $1.94$	& $166.08$	& $68.0$	& $1.84$	& $95.01$\\
		2N2D4K	& $113.2$	& $0.89$	& $431.02$	& $64.4$	& $0.8$	& $348.56$\\
		2N3D3K	& $120.0$	& $1.8$	& $223.69$	& $71.6$	& $1.6$	& $167.72$\\
		2N3D4K	& $115.6$	& $0.27$	& $782.18$	& $13.6$	& $0.1$	& $93.92$\\
		2N4D3K	& $91.6$	& $0.31$	& $184.21$	& $36.2$	& $0.16$	& $108.0$\\
		2N4D4K	& $90.4$	& $0.25$	& $965.26$	& $8.6$	& $0.08$	& $93.7$\\
		\hline
		3N2D3K	& $138.4$	& $7.01$	& $404.59$	& $123.2$	& $6.14$	& $424.47$\\
		3N2D4K	& $150.0$	& $4.7$	& $879.48$	& $62.8$	& $3.99$	& $751.33$\\
		3N3D3K	& $150.0$	& $2.33$	& $301.63$	& $67.0$	& $1.76$	& $160.05$\\
		3N3D4K	& $150.0$	& $0.82$	& $1469.97$	& $66.6$	& $0.35$	& $906.26$\\
		3N4D3K	& $150.0$	& $2.07$	& $354.48$	& $37.2$	& $1.63$	& $110.82$\\
		3N4D4K	& $150.0$	& $1.06$	& $3295.09$	& $37.8$	& $0.56$	& $820.42$\\
		\hline
		4N2D3K	& $150.0$	& $5.01$	& $311.78$	& $103.2$	& $3.66$	& $262.33$\\
		4N2D4K	& $150.0$	& $7.58$	& $1319.14$	& $93.8$	& $5.71$	& $1346.59$\\
		4N3D3K	& $150.0$	& $1.32$	& $317.69$	& $6.6$	& $0.15$	& $16.75$\\
		\multicolumn{1}{c}{}&\multicolumn{3}{|c}{SG} & \multicolumn{3}{|c}{BTM} \\  
		Clustering & $ \overline{t} $ &   $ \overline{\text{rel. DG}} $  & $\overline{T_\text{comp}}$ & $ \overline{t} $ &   $ \overline{\text{rel. DG}} $  & $\overline{T_\text{comp}}$ \\
		\hline
		\textbf{Mean}	& $\vec{54.48}$	& $\vec{1.81}$	& $\vec{483.22}$&  &   & \\
		\hline
		4N3D4K	& $150.0$	& $1.55$	& $2786.47$	& $65.8$	& $0.53$	& $2046.2$\\
		4N4D3K	& $150.0$	& $1.57$	& $441.69$	& $35.0$	& $0.42$	& $122.26$\\
		4N4D4K	& $150.0$	& $1.7$	& $2593.41$	& $7.2$	& $0.14$	& $118.24$\\		
		\hline\hline
		\multicolumn{1}{c}{}&\multicolumn{3}{|c|}{QNDA} &  &   & \\  
		Clustering & $ \overline{t} $ &   $ \overline{\text{rel. DG}} $  & $\overline{T_\text{comp}}$ &  &   & \\
		\hline
		\textbf{Mean}	& $\vec{54.48}$	& $\vec{1.81}$	& $\vec{483.22}$&  &   & \\
		\hline
		2N2D3K	& $63.2$	& $1.82$	& $92.57$&  &   & \\
		2N2D4K	& $62.2$	& $0.73$	& $345.59$&  &   & \\
		2N3D3K	& $62.2$	& $1.58$	& $150.11$&  &   & \\
		2N3D4K	& $5.0$	& $0.06$	& $34.76$&  &   & \\
		2N4D3K	& $32.8$	& $0.12$	& $97.89$&  &   & \\
		2N4D4K	& $4.8$	& $0.08$	& $47.19$&  &   & \\
		\hline
		3N2D3K	& $121.0$	& $6.21$	& $412.26$&  &   & \\
		3N2D4K	& $64.2$	& $3.98$	& $747.3$&  &   & \\
		3N3D3K	& $64.0$	& $1.71$	& $151.79$&  &   & \\
		3N3D4K	& $63.8$	& $0.28$	& $858.76$&  &   & \\
		3N4D3K	& $35.0$	& $1.36$	& $111.44$&  &   & \\
		3N4D4K	& $37.2$	& $0.46$	& $731.33$&  &   & \\
		\hline
		4N2D3K	& $94.6$	& $3.61$	& $249.28$&  &   & \\
		4N2D4K	& $93.8$	& $5.68$	& $1281.33$&  &   & \\
		4N3D3K	& $9.4$	& $0.17$	& $22.76$&  &   & \\
		4N3D4K	& $66.0$	& $0.49$	& $1901.69$&  &   & \\
		4N4D3K	& $37.4$	& $0.41$	& $129.74$&  &   & \\
		4N4D4K	& $9.8$	& $0.17$	& $178.24$&  &   & \\				
		\hline								
	\end{tabular}
 \end{table*}
\vfill

\end{document}